\documentclass[11pt]{lmcs}
\pdfoutput=1

\usepackage{lastpage}
\lmcsdoi{20}{4}{7}
\lmcsheading{}{\pageref{LastPage}}{}{}%
{Nov.~23,~2022}{Oct.~10,~2024}{}

\usepackage[utf8]{inputenc}
\usepackage{latexsym,amsmath,amssymb,amsthm,stmaryrd,url}
\usepackage{graphicx}
\usepackage[inline]{enumitem}

\newcommand\C{\mathfrak{C}}
\newcommand\Cb{\mathfrak{D}}
\newcommand\Val{\mathbf V}
\newcommand\Lform{{\mathcal L}}
\newcommand\wk{\mathrm{w}}
\newcommand\pw{\mathrm{p}}

\newcommand\cb[1]{\mathbf{#1}} 

\newcommand{\Topcat}{\cb{Top}}
\newcommand{\Contcat}{\cb{Cont}}
\newcommand{\SCcat}{\cb{SComp}}

\newcommand{\real}{\mathbb{R}}
\newcommand\Rp{\real_+}
\newcommand\creal{\overline{\real}_+}

\newcommand\limp{\mathrel{\Rightarrow}}
\newcommand\dsup{\sup\nolimits^{\scriptstyle\uparrow}}

\newcommand\eqdef{\mathrel{\buildrel \text{def}\over=}}
\newcommand\upc{\mathop{\uparrow}}
\newcommand\dc{\mathop{\downarrow}}
\newcommand\uuarrow{\rlap{$\uparrow$}\raise.5ex\hbox{$\uparrow$}}
\newcommand\ddarrow{\rlap{$\downarrow$}\raise.5ex\hbox{$\downarrow$}}
\newcommand\Sm{{\mathcal S}}
\newcommand\Min{\mathop{\text{Min}}}
\newcommand\diff{\smallsetminus}
\newcommand{\interior}[1]{int ({#1})} 

\newcommand{\blue}[1]{\textcolor{black}{#1}}

\newcommand\dee{\;d} 
\newcommand{\BLUE}[1]{\textcolor{black}{#1}}


\begin{document}

\title[A cone-theoretic barycenter existence theorem]{A cone-theoretic barycenter existence theorem}

\author[J.~Goubault-Larrecq]{Jean Goubault-Larrecq \lmcsorcid{0000-0001-5879-3304}}[a]
\author[X.~Jia]{Xiaodong Jia\lmcsorcid{0000-0001-9310-6143}}[b]

\address{Universit\'e Paris-Saclay, CNRS, ENS Paris-Saclay, Laboratoire M\'ethodes Formelles, 91190, Gif-sur-Yvette, France}	
\email{jgl@lmf.cnrs.fr} 

\address{School of Mathematics, Hunan University, Changsha, Hunan, 410082, China}
\email{jiaxiaodong@hnu.edu.cn (Corresponding author)}


\begin{abstract}
  We show that every continuous valuation on a locally convex, locally
  convex-compact, sober topological cone $\C$ has a barycenter.  This
  barycenter is unique, and the barycenter map $\beta$ is continuous,
  hence is the structure map of a $\Val_\wk$-algebra, i.e., an
  Eilenberg-Moore algebra of the extended valuation monad on the
  category of $T_0$ topological spaces; it is, in fact, the unique
  $\Val_\wk$-algebra that induces the cone structure on~$\C$.
\end{abstract}

\maketitle


\section{Introduction}
\label{sec:introduction}

It is well-known that every probability measure $\mu$ defined on a
compact convex subset of a Hausdorff topological vector space has a
unique barycenter $\beta (\mu)$, and that the barycenter map $\beta$
is continuous \cite[Chapter~6, 26.3]{Choquet:analysis:2}.  Our purpose
is to establish a similar theorem in the setting of topological cones,
in the sense of Reinhold Heckmann \cite{heckmann96} and Klaus Keimel
\cite{Keimel:topcones2}; those are definitely not limited to cones
that one can embed in vector spaces.

\paragraph{Applications.}
Such a barycenter existence theorem is strongly linked to the question
of the elucidation of the algebras of the monad $\Val_{\wk}$ of
continuous valuations on the category $\Topcat_0$ of $T_0$ topological
spaces, as well as of the monads of probability and subprobability
valuations.  (See \cite{CEK:prob:scomp,GLJ:bary}.  We will introduce
$\Val_\wk$ below.)  As a result, we will be able to show that a rich
collection of topological cones give rise to such algebras.

Investigating the structure of algebras of monads is not just an
interesting, and important, mathematical question, but also has
application in semantics.  For example, knowing that the category of
algebras of the subprobability valuation monad on the category of
continuous dcpos is isomorphic to the category of so-called continuous
Kegelspitzen is crucial to give an adequate semantics to the
variational quantum programming language of \cite{10.1145/3498687},
where barycenter maps are used to link the category of some special
Von-Neumann algebras to the Kleisli category of the so-called minimal
valuations monad.

\paragraph{Related work.}
Ben Cohen, Mart\'{\i}n Escard{\'o} and Klaus Keimel
\cite{CEK:prob:scomp} first asked the question of characterizing the
algebras of $\Val_\wk$ on $\Topcat_0$, and then on the subcategory
$\SCcat$ of stably compact spaces and continuous maps.  They also
observed how the question was deeply linked to a notion of barycenter
adapted from Choquet's work.  There seems to be a gap in their
Proposition~2, however (see Proposition~\ref{prop:Valg} below), and
\cite{GLJ:bary} is probably now a better reference.  Both papers only
give \emph{examples} of spaces where such barycenters exist, as well
as of algebras of $\Val_\wk$, but no general existence theorem.

Keimel also fully characterized the algebras of the monad of
probability measures over compact ordered spaces
\cite{Keimel:Valg:scomp}.  This can be seen as a close result, because
of the tight relationship between compact ordered spaces and stably
compact spaces.  But the category of compact ordered spaces and
continuous maps is equivalent to the category of stably compact spaces
and \emph{perfect} maps, not continuous maps.  Keimel's result itself
extends \'Swirszcz's characterization of the algebras of the monad of
probability measures on compact Hausdorff spaces
\cite{Swirszcz:convex}.

Returning to categories of spaces that are not necessarily Hausdorff,
with continuous, not necessarily perfect maps, we do have some results
for certain submonads of $\Val_\wk$: this includes the monads of
so-called simple valuations and of point-continuous valuations
\cite{GLJ:bary,heckmann96}, for which we have general barycenter
existence theorems, and complete characterizations of the
corresponding algebras.  Since point-continuous valuations and
continuous valuations coincide on continuous dcpos, those results
encompass the case of continuous Kegelspitzen mentioned earlier.

Our purpose here is to give a barycenter existence theorem that would
need assumptions that are as weak as possible.  We initially looked
for variants of arguments due to Choquet \cite{Choquet:repr}, Edwards
\cite{Edwards:repr}, and Roberts \cite{Roberts:embedding}, among
others.  These arguments work as follows.  Define a partial ordering
between probability distributions by $\mu_1 \preceq \mu_2$ if and only
if $\int h \dee\mu_1 \geq \int h \dee\mu_2$ for every convex map $h$ in some
class (typically, bounded, continuous maps); this formalizes the idea
that $\mu_2$ has a more concentrated distribution of mass than
$\mu_1$, while keeping the same center of mass.  Then, build a
continuous map $F$ such that $\mu \preceq F (\mu)$ for all probability
distributions $\mu$, and such that $F (\mu)=\mu$ if and only if $\mu$
is a Dirac measure $\delta_x$.  Finally, using a compactness argument,
extract a convergent subsequence from the sequence
${(F^n (\mu))}_{n \in \mathbb N}$: this must converge to some Dirac
measure $\delta_x$, and by construction $x$ will be a barycenter of
$\mu$.  Details may vary.  While there are immediate obstructions to
generalizing that form of argument to a non-Hausdorff setting (e.g.,
limits are not unique), it is possible to do so, but the result we
obtained by following this path was disappointing: a long, complicated
proof, which required very strong assumptions in the end.

In comparison, we will use a much simpler strategy: we reduce the
question to the already known construction of barycenters of
point-continuous valuations, by embedding our cone into a larger one,
the upper powercone introduced by Keimel
\cite[Section~11]{Keimel:topcones2}.  We only need two additional
lemmata: one to actually show that there is such a cone embedding,
Lemma~\ref{lemma:Lambda:bar:Qcvx}, and one to transfer any barycenter
in the larger cone back to the smaller cone, Lemma~\ref{lemma:jia}.

\section{Cones, continuous valuations, barycenters}
\label{sec:cones-cont-valu}

We refer the reader to \cite{GHKLMS:contlatt,JGL-topology} for more
information on domain theory and topology, especially non-Hausdorff
topology.

A \emph{dcpo} (short for directed-complete partial order) is a poset
$(X, \leq)$ such that every directed family $D$ has a supremum
$\sup D$; we write $\dsup D$ to stress the fact that $D$ is directed.
\blue{ A subset $A$ of $D$ is \emph{upwards-closed} if and only if
  $x \in A$ and $x \leq y$ imply $y \in A$.  We will also write
  $\upc A$ for the upward closure
  $\{y \in X \mid x \leq y \text{ for some }x \in A\}$, and $\upc x$
  instead of $\upc \{x\}$,  \blue{and similarly for the downward closure
    $\dc A$, and $\dc x$.}  The \emph{Scott topology} on $(X, \leq)$
  consists of Scott-open subsets of $X$, which are defined to be
  upwards-closed and inaccessible by directed suprema. That is, an
  upwards-closed subset $U$ is \emph{Scott-open} if and only if every
  directed family $D$ intersects $U$ if $\dsup D \in U$.  A
  Scott-continuous map between dcpos is a monotonic map that preserves
  directed suprema; equivalently, a continuous map, provided the given
  dcpos are given the Scott topology. }

We write $\Rp$ for the set of non-negative real numbers, and $\creal$
for $\Rp \cup \{\infty\}$, where $\infty$ is larger than any element
of $\Rp$. \BLUE{$\creal$ is a dcpo in the natural ordering $\leq$ on real numbers, 
$\Rp$ is not.  We equip the dcpo~$\creal$ with the Scott topology;} its
non-empty open subsets are the half-open intervals $]a, \infty]$,
$a \in \Rp$.  (And we give $\Rp$ the subspace topology.)  With that
topology, the continuous maps from a topological space to $\creal$ are
usually called the \emph{lower semicontinuous maps}.

The way-below relation $\ll$ on a dcpo $X$ is defined by $x \ll y$ if
and only if for every directed family $D$ such that $y \leq \dsup D$,
there is a $z \in D$ such that $x \leq z$.  We write $\ddarrow y$ for
$\{x \in X \mid x \ll y\}$.  A \emph{continuous dcpo} is a dcpo in
which $\ddarrow y$ is directed and $\dsup \ddarrow y = y$ for every
point $y$.

\blue{The \emph{specialization preordering} $\leq$ of a topological
  space $X$ is defined by $x \leq y$ if and only if every open
  neighborhood of $x$ contains $y$.  The space $X$ is $T_0$ if and
  only if $\leq$ is antisymmetric.  The \emph{saturation} of a subset
  $A$ of $X$ is the intersection of its open neighborhoods, and
  coincides with $\upc A$.  The closure of a singleton $\{x\}$ is
  equal to $\dc x$, where $\dc$ is relative to the specialization
  preordering.  Any continuous map between topological spaces is
  monotonic with respect to the underlying specialization preorderings.}

\blue{A subset $K$ of a topological space $X$ is \emph{compact} if and
  only if every open cover of $K$ contains a finite subcover; no
  separation axiom is intended here.  If $K$ is compact, then $\upc K$
  is both compact and saturated.  A topological space is \emph{locally
    compact} if every point in it has a neighborhood base of compact
  subsets, namely if for every point $x$ and every open neighborhood
  $U$ of $x$, there is a compact set $Q$ such that
  $x \in \interior Q \subseteq Q \subseteq U$.  (We write
  $\interior K$ for the \emph{interior} of $K$.)  We may equivalently
  require $Q$ to be compact \emph{saturated}, by replacing $Q$ by
  $\upc Q$.  A space is \emph{sober} if every nonempty irreducible
  closed subset in it is equal to $\dc x$ for some unique point $x$. }
Every continuous dcpo is locally compact and sober in its Scott
topology~\blue{\cite[Corollary II-1.13]{GHKLMS:contlatt}.}


Following Heckmann \cite{heckmann96} and Keimel
\cite{Keimel:topcones2}, a \emph{cone} $\C$ is an additive commutative
monoid with an action $a, x \mapsto a \cdot x$ of the semi-ring $\Rp$
on $\C$ (\emph{scalar multiplication}).  The latter means that the
following laws are satisfied:
\[
  \begin{array}{ccc}
    0 \cdot x = 0 & (ab) \cdot x = a \cdot (b \cdot x) & 1 \cdot x=x\\
    a \cdot 0=0 & a \cdot (x+y) = a \cdot x+a \cdot y & (a+b) \cdot x = a \cdot x+b \cdot x.
  \end{array}
\]
We will often write $ax$ instead of $a \cdot x$.  In other words, a
cone satisfies the same axioms as a vector space, barring those
mentioning subtraction.  $\Rp$ is a cone, and embeds (\blue{as a
  cone}) into the real vector space $\real$.  $\creal$ is a cone, but
embeds in no real vector space, because it fails the cancellation law
$x+y=x+z \limp y=z$.  None of the other cones we will consider embed
in any real vector space, for the same reason.

\BLUE{In this paper, a \emph{topological cone} is a cone with a $T_0$ topology }
that makes both addition and scalar multiplication jointly continuous---where we
recall that the topology on $\Rp$ is induced by the Scott topology on
$\creal$, which is not Hausdorff.  A \emph{semitopological cone} is
defined similarly, except that addition is only required to be
separately continuous.  \blue{Since continuous maps are monotonic with
  respect to the specialization preordering $\leq$}, in a
semitopological cone, $ax \leq bx$ whenever $a \leq b$; in particular,
$0$ is the least element.  Hence a non-trivial semitopological cone
\blue{$\C$} is never $T_1$, but is always compact.  \blue{Compactness
  comes from the fact that every open cover ${(U_i)}_{i \in I}$ will
  be such that $0 \in U_i$ for some $i \in I$, and then $U_i = \C$,
  since $U_i$ is upwards-closed and $0$ is the least element of $\C$.}

A \emph{d-cone} is a cone with a partial ordering that turns it into a
dcpo, and such that addition and scalar multiplication are both
Scott-continuous.  \blue{A continuous d-cone is a d-cone that is also
  a continuous dcpo.}  Every continuous d-cone is a topological cone
in the Scott topology.  For example, $\creal$ is a continuous d-cone.

A subset $A$ of a cone $\C$ is \emph{convex} if and only if
$ax+(1-a)y \in A$ for all $x, y \in A$ and $a \in [0, 1]$.  A
semitopological cone is \emph{weakly locally convex}, resp.\
\emph{locally convex}, resp.\ \emph{locally convex-compact}, if and
only if every point has a base of convex, resp.\ convex open, resp.\
convex compact, neighborhoods.  Every continuous d-cone is not only
topological, but also locally convex and locally convex-compact
\cite[Lemma~6.12]{Keimel:topcones2}.

  \begin{exa}
    \label{example-semilattice-cones}
    Here are a few examples to bear in mind.  We will give a few
    others below.  The point of \ref{exa:cone:weird} and
    \ref{exa:cone:weird:upper} is twofold: to show that cones can be
    very different from real vector spaces; and to show that there are
    sober, locally convex, locally convex-compact topological cones
    that are not continuous d-cones.
    \begin{enumerate}[series=exacone,label=(\arabic*)]
    \item For every space $X$, let us consider the dcpo $\Lform X$ of
      all lower semicontinuous maps from $X$ to $\creal$, namely,
      those continuous maps from $X$ to $\creal$ equipped with the
      Scott topology.  With the obvious addition and scalar
      multiplication, and the Scott topology of the pointwise
      ordering, $\Lform X$ is a semitopological cone.

      If $X$ is core-compact, namely, if its lattice of open sets is a
      continuous dcpo (every locally compact space is core-compact,
      and the two notions are equivalent for sober spaces), then
      $\Lform X$ is a continuous d-cone \cite[Example~3.5]{GLJ:bary}.
      In particular, it is topological, sober, locally convex, and
      locally convex-compact.
    \item\label{exa:cone:weird} This one is due to Heckmann
      \cite[Section~6.1]{heckmann96}, see also
      \cite[Example~3.6]{GLJ:bary}.  Every sup-semilattice $L$ with a
      least element gives rise to a semitopological cone in its Scott
      topology by defining $x+y$ as $x \vee y$, $0$ as the least
      element $\bot$, and $a \cdot x$ as $\bot$ if $a=0$, and as $x$
      otherwise; it is a continuous d-cone (hence topological, locally
      convex, and locally convex-compact) if $L$ is a continuous
      complete lattice.

      For general $L$, it remains that every upwards-closed subset is
      convex. Hence each such sup-semilattice $L$ is a locally convex
      semitopological cone in its Scott topology. It is locally
      convex-compact if $L$ is also locally compact in the Scott
      topology.  In that case, $L$ is also a topological cone. This is
      because on a core-compact poset $L$, the Scott topology on
      $L\times L$ is the same as the product topology of the Scott
      topology on $L$~\cite[Theorem II-4.13]{GHKLMS:contlatt}.

      \noindent
      \begin{minipage}{0.85\linewidth}
        In particular, consider any locally compact sup-semilattice
        $L$ with bottom that is not a continuous dcpo.  Then $L$ is
        topological, locally convex, and locally convex-compact, but
        not a continuous d-cone.  The simplest example is given by
        putting two copies of $\mathbb N$ side by side, with extra
        bottom and top elements, as shown on the right.  Its nonempty
        proper Scott open subsets are those of the form
        $\upc \{m, n'\}$ where $m, n \in \mathbb N$, and they are all
        compact, so $L$ is locally compact.  One can even check that
        it is a quasi-continuous dcpo
        \cite[Definition~III-3.2]{GHKLMS:contlatt}, which not only
        implies that it is locally compact, but also sober
        \cite[Proposition~III-3.7]{GHKLMS:contlatt}.  However, $L$ is
        not a continuous dcpo, since the way-below relation is
        characterized by $x \ll y$ if and only if $x=\bot$.
      \end{minipage}
      \begin{minipage}{0.14\linewidth}
        \includegraphics[scale=0.2]{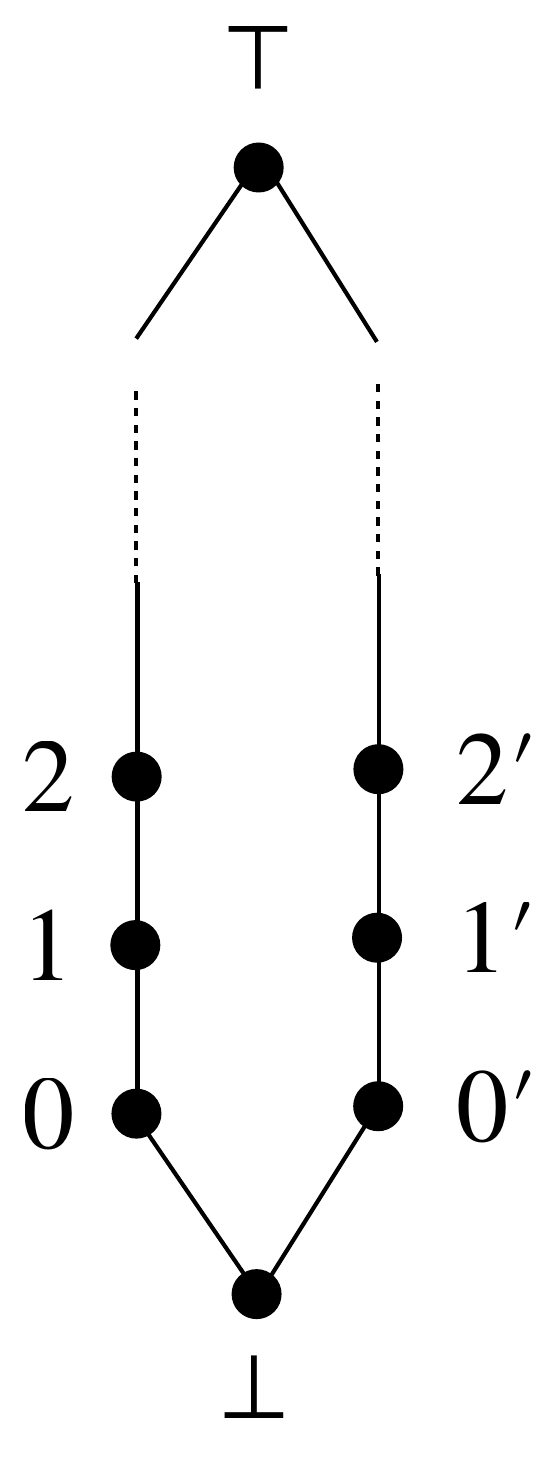}
      \end{minipage}
    \item\label{exa:cone:weird:upper} Here is a variant on the
      previous example.  We consider a complete lattice $L$ with its
      \emph{upper topology}, which is the coarsest topology that makes
      the sets $\dc x$ closed.  The specialization ordering is the
      ordering $\leq$ of $L$.  By a result of Schalk
      \cite[Proposition~1.7]{Schalk:PhD}, $L$ is sober.  We equip $L$
      with the same cone structure as in~\ref{exa:cone:weird}.  Then
      $L$ is topological, since
      ${(\_ + \_)}^{-1} (\dc x) = \dc x \times \dc x$ is closed in the
      product topology for every $x \in L$, and similarly for
      ${(\_ \cdot \_)}^{-1} (\dc x) = (\Rp\times \dc x) \cup (\{0\}
      \times L)$.  As in~\ref{exa:cone:weird}, every upwards-closed
      subset of $L$ is convex, so $L$ is locally convex, and $L$ is
      locally convex-compact if and only if it is locally compact.

      In particular, the lattice $\Gamma (X)$ of closed subsets of any
      topological space $X$, with the upper topology of inclusion
      (also known as the \emph{lower Vietoris topology}), union as
      addition, and empty set as zero, is a sober, locally convex
      topological cone.  It is locally compact, hence locally
      convex-compact, if $X$ is itself locally compact
      \cite[Proposition~6.11]{Schalk:PhD}.  But, say,
      $\Gamma (\real)$, where $\real$ is taken with its usual metric
      topology (a locally compact space), is not a continuous dcpo (or
      a quasi-continuous dcpo) with its Scott topology; the point is
      that the upper topology differs from the Scott topology: by a
      result of Chen, Kou and Lyu
      \cite[Proposition~3.15]{CKL:closed:topol}, the only $T_1$ spaces
      $X$ such that the upper and Scott topologies coincide on
      $\Gamma (X)$ are discrete.

    \end{enumerate}
\end{exa}
\noindent 
We write $\C^*$ for the cone of linear lower semicontinuous maps
$\Lambda \colon \C \to \creal$, where linear means preserving sums and
scalar products.  Every locally convex $T_0$ semitopological cone is
\emph{linearly separated}, in the sense that for every pair of points
$x$, $y$ such that $x \not\geq y$, there is a map $\Lambda \in \C^*$
such that $\Lambda (x) < \Lambda (y)$
\cite[Corollary~9.3]{Keimel:topcones2}.  By multiplying by an
appropriate constant, we may require $\Lambda (x) \leq 1$ and
$\Lambda (y) > 1$.

A \emph{continuous valuation} $\nu$ on a space $X$ is a map from the
lattice of open subsets of $X$ to $\creal$ such that
$\nu (\emptyset)=0$,
$\nu (U \cup V) + \nu (U \cap V) = \nu (U) + \nu (V)$ for all open
sets $U$ and $V$, and such that $\nu$ is Scott-continuous
\cite{Jones:proba}.  This is a close cousin of the notion of measure:
every continuous valuation extends to a $\tau$-smooth measure
on the Borel $\sigma$-algebra of $X$ provided that $X$ is LCS-complete
\cite[Theorem~1]{dBGLJL:LCS}, and that includes the case of all
locally compact sober spaces.  Additionally, there is a notion of
integral of maps $h \in \Lform X$ with respect to continuous
valuations $\nu$ on $X$, with the usual properties.  We will use the
following four:
\begin{itemize}
\item $\int_{x \in X} h (x) \dee\delta_{x_0} = h (x_0)$, where
  $\delta_{x_0}$ is the \emph{Dirac valuation} at $x_0$, defined by
  $\delta_{x_0} (U) \eqdef 1$ if $x_0 \in U$, $0$ otherwise;
\blue{\item  $\int_{x \in X} \chi_U (x) \dee \nu  = \nu(U)$, where $\chi_U$
is the characteristic function of an open set $U$, defined as 
$\chi_U(x) = 1$ if $x \in U$, $0$ otherwise;}
\item the \emph{change of variable} formula: for every continuous map
  $f \colon X \to Y$,
  $\int_{y \in Y} h (y) \dee f[\nu] = \int_{x \in X} h (f (x)) \dee\nu$; here
  $f [\nu]$ is the \emph{image valuation} defined by
  $f [\nu] (V) \eqdef \nu (f^{-1} (V))$;
\item integration is linear in the integrated function $h$.
\end{itemize}
\noindent 
We write $\Val X$ for the set of all continuous valuations on $X$.
The \emph{weak topology} on $\Val X$ is the coarsest one that makes
$\nu \mapsto \nu (U)$ continuous for every open set $U$.  In other
words, a subbase of the weak topology is given by the sets
$[U > r] \eqdef \{\nu \in \Val X \mid \nu (U) > r\}$, where $U$ ranges
over the open subsets of $X$ and $r \in \Rp$.  We write $\Val_\wk X$
for $\Val X$ with the weak topology.  The specialization ordering on
$\Val_\wk X$ is the so-called \emph{stochastic ordering}: for all
continuous valuations $\mu$ and $\nu$, $\mu \leq \nu$ if and only if
$\mu(U) \leq \nu(U)$ for every open subset $U$ of $X$.  $\Val_\wk$
extends to an endofunctor on the category $\Topcat_0$ of $T_0$
topological spaces, and its action on morphisms $f \colon X \to Y$ is
given by $\Val_\wk f (\nu) \eqdef f [\nu]$.

\blue{%
  \begin{exa}
    \label{exa:cone:Val}
    Here are few more examples of semitopological cones.
    \begin{enumerate}[resume*=exacone]
    \item\label{exa:V} With the obvious addition and scalar
      multiplication, $\Val_\wk X$ is a locally convex, sober,
      topological cone \cite[Proposition~3.13]{GLJ:bary}.  This is
      really a special case of the next item, considering that
      $\Val_\wk X$ is isomorphic to the dual cone ${(\Lform X)}^*$
      (see next item) as semitopological cones.
    \item\label{exa:C*} We equip $\C^*$ with the \emph{weak$^*$upper
        topology}, namely the coarsest one that makes the functions
      $\text{ev}_x \colon \Lambda \mapsto \Lambda (x)$ lower
      semicontinuous from $\C^*$ to $\creal$, for each $x \in \C$;
      $\creal$ again has the Scott topology of its ordering.  (The
      weak$^*$upper topology is also known as the \emph{weak$^*$-Scott
        topology} \cite{Plotkin:alaoglu}.)  The cone $\C^*$ is the
      \emph{dual cone} of $\C$ \cite{Keimel:topcones2}; see also
      \cite[Section~3]{CEK:prob:scomp} or Example~3.7~(ii) of
      \cite{GLJ:bary}.  A dual cone is always topological
      \cite[Example~3.7~(ii)]{GLJ:bary}.  It is also locally convex,
      since it has a base of open sets of the form
      $\{\Lambda \in \C^* \mid \Lambda (x_1) > 1, \cdots, \Lambda
      (x_n) > 1\}$, where $x_1$, \ldots, $x_n$ ranges over the finite
      lists of points of $\C$, and such sets are convex.  Finally,
      $\C^*$ is sober, using the following argument, adapted from an
      argument of Heckmann's \cite[Proposition~5.1]{heckmann96}.  For
      every sober space $Y$, the space $[X \to Y]_{\mathrm{p}}$ of \BLUE{all
      continuous maps} from $X$ to $Y$, with the pointwise
      topology, is sober by \cite[Lemma~5.8]{Tix:bewertung}; the
      pointwise topology is the coarsest that makes $h \mapsto h (x)$
      continuous for every $x \in X$.  By letting $Y \eqdef \creal$
      with the Scott topology, $[X \to \creal]_{\mathrm{p}}$ is sober.
      Then $\C^*$ occurs as a $T_0$-equalizer subspace of
      $[X \to \creal]_{\mathrm{p}}$, namely as the equalizer of the
      continuous maps
      $h \mapsto ({(h (ax))}_{a \in \Rp, x \in \C}, {(h (x+y))}_{x, y
        \in \C})$ and
      $h \mapsto ({(a h (x))}_{a \in \Rp, x \in \C}, {(h
        (x)+h(y))}_{x, y \in \C})$ with codomain the $T_0$ space
      $\creal^{\Rp \times \C} \times \creal^{\C^2}$.  And any
      $T_0$-equalizer subspace of a sober space is sober
      \cite[Lemma~5.9]{Tix:bewertung}.
    \end{enumerate}
  \end{exa}}
\noindent 
$\Val_\wk$ is the functor part of a monad on $\Topcat_0$, whose unit
$\eta_X$ maps every point $x \in X$ to $\delta_x$, \blue{and the
  multiplication map at $X$ sends a continuous valuation
  $\varphi \in \Val_\wk \Val_\wk X$ to
  $(U\mapsto \int_{\nu\in \Val_\wk X} \nu(U) d \varphi)$.}  We are
interested in its algebras $\alpha \colon \Val_\wk X \to X$.  By
Lemma~4.6 of \cite{GLJ:bary}, any such algebra induces a cone
structure on $X$ by $x+y \eqdef \alpha (\delta_x+\delta_y)$,
$a \cdot x \eqdef \alpha (a \delta_x)$, which turns $X$ into a weakly
locally convex, sober, topological cone
\cite[Proposition~4.9]{GLJ:bary}.  Additionally, with this induced
cone structure, $\alpha$ is a linear map, and for every
$\nu \in \Val_\wk X$, $\alpha (\nu)$ is a barycenter of $\nu$.

Following \cite{CEK:prob:scomp,GLJ:bary}, we define a
\emph{barycenter} of $\nu \in \Val_\wk \C$, where $\C$ is a
semitopological cone, as any point $x_0$ of $\C$ such that
$\Lambda (x_0) = \int_{x \in \C} \Lambda (x) \dee\nu$ for every
$\Lambda \in \C^*$.  If $\C$ is linearly separated and $T_0$, then
barycenters are unique if they exist.  This definition is similar to
Choquet's classic definition \cite[Chapter~6,
26.2]{Choquet:analysis:2}.  The main difference is that our linear
maps are only required to be lower semicontinuous; also, we have
replaced vector spaces by cones, and the definition applies to all
continuous valuations, not just probability valuations (namely,
valuations $\nu$ such that $\nu (\C)=1$).

\blue{The following example should explain how barycenters generalize
  the usual notion of barycenters of points $x_i$ with weights $a_i$.
\begin{exa}
  \label{exa:bary:simple}
  A \emph{simple valuation} is any finite linear combination
  $\nu \eqdef \sum_{i=1}^n a_i \delta_{x_i}$ of Dirac valuations, with
  $a_i \in \Rp$.  If all the points $x_i$ are taken from a
  semitopological cone $\C$, then $\sum_{i=1}^n a_i \cdot x_i$ is a
  barycenter of $\nu$, by the linearity of integration and of the maps
  $\Lambda$ involved in the definition \cite[Example~4.3]{GLJ:bary}.
\end{exa}
\begin{exa}
  \label{exa:bary}
  Here are additional, more sophisticated examples.
  \begin{enumerate}[label=(\arabic*)]
  \item For every core-compact space $X$, every continuous valuation
    on $\Lform X$ has a barycenter.  This is really a consequence of
    Proposition~\ref{prop:cont:Valg} below, plus the fact that
    $\Lform X$ is a continuous d-cone in that case.  Explicitly, one
    can check that the map
    $x \mapsto \int_{h \in \Lform X} h (x) \dee\nu$ is that barycenter
    \cite[Proposition~4.15]{GLJ:bary}.
  \item\label{exa:supp:weird} If $L$ is a complete lattice with its
    Scott topology (see
    Example~\ref{example-semilattice-cones}~\ref{exa:cone:weird}),
    then every continuous valuation $\mu$ on $L$ has a barycenter,
    which is the supremum of the elements of its support
    $\mathrm{supp}\;\mu$ \cite[Example~4.4]{GLJ:bary}; the
    \emph{support} $\mathrm{supp}\; \mu$ of $\mu$ is defined as the
    complement of the largest open set $U$ such that $\mu(U) = 0$.
  \item\label{exa:supp:weird:upper} Let $L$ be a complete lattice with
    its upper topology (see
    Example~\ref{example-semilattice-cones}~\ref{exa:cone:weird:upper}).
    Just like in the previous item, every continuous valuation $\mu$
    on $L$ has a barycenter, which is the supremum of
    $\mathrm{supp}\;\mu$.  The argument follows the same lines as
    \cite[Examples~3.15, 3.20, 4.4]{GLJ:bary}.  Each non-empty convex
    subset of $L$ is directed.  Every non-empty closed
    convex subset $C$ must then be of the form $\dc x$ for some
    $x \in L$: namely, $x = \dsup C$, which must be in $C$ since $C$
    is closed, hence Scott-closed (note that the upper topology is
    coarser than the Scott topology).  The elements $\Lambda$ of $L^*$
    must be such that $\Lambda^{-1} ([0, 1])$ is closed, non-empty,
    and convex, hence of the form $\dc x_0$ for some $x_0 \in L$.  For
    every $x \in L$, the fact that
    $\Lambda (x) = \Lambda (r \cdot x) = r \Lambda (x)$ for every
    $r > 0$ implies that $\Lambda (x)$ can only be equal to $0$ or to
    $\infty$.  Hence $\Lambda = \infty \cdot \chi_{L \diff \dc x_0}$.
    Conversely, any map of this form is in $L^*$.  Using the
    definition of barycenters, a point $x \in L$ is a barycenter of
    $\mu$ if and only if for every $x_0 \in L$, \BLUE{the inequality $x \leq x_0$ is equivalent
    to that $\infty \cdot \mu (L \diff \dc x_0)=0$, to that
    $\mathrm{supp}\;\mu \subseteq \dc x_0$, and to that 
    $\sup (\mathrm{supp}\;\mu) \leq x_0$, whence the conclusion.}
  \item\label{exa:C*:bary} The barycenters of continuous
    valuations $\nu$ on dual cones $\C^*$ do not exist in general. However, if
    $\Lambda_0$ is a barycenter of a continuous valuation $\nu$ on
    $\C^*$, then
    $\varphi (\Lambda_0) = \int_{\Lambda \in \C^*} \varphi (\Lambda)
    d\nu$ for every $\varphi \in \C^{**}$ by definition of
    barycenters.  For each point $x \in \C$, the map
    $\varphi \colon \Lambda \mapsto \Lambda (x)$ is an element of
    $\C^{**}$, so
    $\Lambda_0 (x) = \int_{\Lambda \in \C^*} \Lambda (x) \dee\nu$ for
    every $x \in \C$.  In other words, if $\nu$ has a barycenter, then
    it must be the map
    $x \mapsto \int_{\Lambda \in \C^*} \Lambda (x) \dee\nu$.  However,
    it is unknown whether that map is a barycenter of $\nu$. Also, that map 
    may fail to be an element of $\C^*$ at all, as it is
    a map that may well fail to be lower semicontinuous.  We will
    return to this in
    Example~\ref{exa:thm:appl}~\ref{exa:C*:bary:yes}.
  \item\label{exa:bary:V} Every continuous valuation $\varphi$ on
    $\Val_\wk X$, for any space $X$ whatsoever, has a barycenter, and
    it is obtained by applying the multiplication of the $\Val_\wk$
    monad to $\varphi$.  Indeed, by general category theory, that
    multiplication is a $\Val_\wk$-algebra (in fact, the free
    $\Val_\wk$-algebra), and we have seen that every
    $\Val_\wk$-algebra is a barycenter map.
  \end{enumerate}
\end{exa}
}
\noindent 
We have seen that every $\Val_\wk$-algebra $\alpha$ induces a cone
structure on its codomain, and that $\alpha$ maps every continuous
valuation to one of its barycenters.  Conversely, we have:
\begin{prop}[Proposition~4.10 of \cite{GLJ:bary}]
  \label{prop:Valg}
  If every continuous valuation $\nu$ on a semitopological cone $\C$
  has a unique barycenter $\beta (\nu)$, and if $\beta$ is continuous,
  then $\beta$ is a $\Val_\wk$-algebra, and is in fact the unique
  $\Val_\wk$-algebra that induces the given cone structure on $\C$.
\end{prop}
\noindent 
The continuity of $\beta$ is essential, and seems to have been
forgotten in \cite[Proposition~2]{CEK:prob:scomp}.

\section{The barycenter existence theorem}
\label{sec:main-theorem}

We start with a topological cone $\C$.  We use a construction due to
Keimel \cite[Section~11]{Keimel:topcones2}: the \emph{upper powercone}
$\Sm \C$ of $\C$ consists of all non-empty convex, compact saturated
subsets of $\C$.  (A \emph{saturated} set is one that is the
intersection of its open neighborhoods, or equivalently an
upwards-closed set, with respect to the specialization preordering.)
Addition is defined by
$Q_1 \underline+ Q_2 \eqdef \upc \{x_1+x_2 \mid x_1 \in Q_1, x_2 \in
Q_2\}$, scalar multiplication by
$a \underline\cdot Q \eqdef \upc \{a \cdot x \mid x \in Q\}$, and the
zero is $\underline 0 \eqdef \C$.  The cone $\Sm \C$ is given the
\emph{upper Vietoris topology}, which is generated by the basic open
sets $\Box U \eqdef \{Q \in \Sm \C \mid Q \subseteq U\}$, where $U$
ranges over the open subsets of $\C$.  Its specialization ordering is
\emph{reverse} inclusion $\supseteq$, and every open set is Scott-open
(with respect to $\supseteq$) when $\C$ is sober.  If additionally
$\C$ is locally convex and locally convex-compact, then $\Sm \C$ is a
continuous d-cone, and the Scott and upper Vietoris topologies
coincide \cite[Theorem~12.6]{Keimel:topcones2}.

The construction $\Sm$ is the functor part of a monad on the category
of
sober topological cones and linear continuous maps as homomorphisms.
This is a consequence of the results in Sections~11 and 12
of \cite{Keimel:topcones2}.  We will use its unit, which we write as
$\eta^{\Sm}_\C \colon \C \to \Sm \C$, and which maps every $x$ to its
upward closure $\upc x$.

For every $\Lambda \in \C^*$, for every $Q \in \Sm \C$,
we let $\overline\Lambda (Q) \eqdef \min_{x \in Q} \Lambda (x)$.
\begin{lem}
  \label{lemma:Lambda:bar:Qcvx}
  Let $\C$ be a topological cone.  For every $\Lambda \in \C^*$,
  $\overline\Lambda$ is an element of $(\Sm \C)^*$, and
  $\overline\Lambda \circ \eta^{\Sm}_\C = \Lambda$.
\end{lem}
\begin{proof}
  $\overline\Lambda$ is continuous because
  ${\overline\Lambda}^{-1} (]r, \infty]) = \Box {\Lambda^{-1} (]r,
    \infty])}$ for every $r \in \Rp$.  For every $a \in \Rp$, for
  every $Q \in \Sm \C$,
  $\overline\Lambda (a \underline\cdot Q) = \min_{x \in Q} \Lambda (a
  \cdot x)$, since \blue{$\Lambda$ is continuous hence monotonic,} and taking the minimum of a monotonic map over the
  upward-closure of a set is the same as taking the minimum on that
  set.  Hence
  $\overline\Lambda (a \underline\cdot Q) = a \cdot \min_{x \in Q} \Lambda (x) = a
  \overline\Lambda (Q)$.  Similarly, for all $Q_1, Q_2 \in \Sm \C$,
  $\overline\Lambda (Q_1 \underline+ Q_2) = \min_{x_1 \in Q_1, x_2 \in
    Q_2} \Lambda (x_1+x_2) = \overline\Lambda (Q_1) + \overline\Lambda
  (Q_2)$.  Finally, $\overline\Lambda \circ \eta^{\Sm}_\C$ maps every
  point $x$ to $\min_{y \in \upc x} \Lambda (y) = \Lambda (x)$,
  because $\Lambda$ is monotonic.
\end{proof}
\noindent 
There is a subspace $\Val_\pw X$ of $\Val_\wk X$, for every space $X$,
consisting of so-called \emph{point-continuous valuations}, and due to
Heckmann \cite{heckmann96}.  It does not really matter to us what they
are exactly.  We will use the following facts:
\begin{itemize}
\item $\Val_\pw$ defines a monad on $\Topcat_0$, with similarly
  defined units and multiplication, and its algebras are exactly the
  weakly locally convex, sober topological cones
  \cite[Theorem~5.11]{GLJ:bary}.
  \blue{This would solve the question of finding the
    $\Val_\wk$-algebras if all continuous valuations were
    point-continuous, but that is not the case, even for continuous
    valuations on dcpos \cite{goubault-jia-2021}.}
\item $\Val_\wk$ restricts to a monad on the subcategory $\Contcat$ of
  continuous dcpos.  ($\Val_\wk X$ coincides with $\Val X$ with the
  Scott topology of the stochastic ordering, for every continuous dcpo
  $X$, by a result of Kirch \cite{kirch93}, and Jones showed that
  $\Val$ is a monad on $\Contcat$ \cite{Jones:proba}---although,
  technically, Kirch was the first one to deal with continuous
  valuations with values in $\creal$.)
\item On the subcategory of continuous dcpos, the $\Val_\pw$ and
  $\Val_\wk$ monads coincide \cite[Theorem~6.9]{heckmann96}.
\end{itemize}
\noindent 
This gives us the following result for free.  We remember that every
continuous d-cone is locally convex, sober, and topological.
\begin{prop}
  \label{prop:cont:Valg}
  For every continuous d-cone $\Cb$, there is a unique
  $\Val_\wk$-algebra \blue{$\alpha \colon \Val_\wk \Cb \to \Cb$} that induces
  the given cone structure on $\Cb$.
\end{prop}

\noindent 
The final ingredient we will need is Lemma~\ref{lemma:jia} below.
Given any non-empty compact saturated subset $Q$ of a space $X$, we
let $\Min Q$ be its set of minimal elements\blue{---as usual, with
  respect to the specialization preordering $\leq$ of $X$}.  It is
well-known that $Q = \upc \Min Q$, namely that every element of $Q$ is
above some minimal element.  The quick proof proceeds by applying
Zorn's Lemma to $Q$ with the opposite ordering $\geq$, noting that
this is an inductive poset: for every chain ${(x_i)}_{i \in I}$ of
elements of $Q$, ${(\dc x_i)}_{i \in I}$ is a chain of closed subsets
\blue{(since $\dc x_i$ is the closure of $\{x_i\} $)} that intersect
$Q$, and since $Q$ is compact, their intersection must also intersect
$Q$, say at $x$; then $x \leq x_i$ for every $i \in I$.
\begin{lem}
  \label{lemma:jia}
  Let $\C$ be a linearly separated $T_0$ semitopological cone and $Q$ be a
  non-empty compact saturated subset of $\C$.  If the map
  $\varphi \colon \Lambda \mapsto \min_{x \in Q} \Lambda (x)$ is
  linear from $\C^*$ to $\creal$, then $Q = \upc x$ for some
  $x \in \C$.
\end{lem}
\begin{proof}
  Let us assume that $Q$ cannot be written as $\upc x$ for any
  $x \in \C$, or equivalently, that $\Min Q$ is not reduced to just
  one point.  (Also, $\Min Q$ is non-empty, since $Q$ is non-empty.)
  For every $x \in \Min Q$, $\Min Q$ contains a distinct, hence
  incomparable, point $\overline x$.

  Since $\C$ is linearly separated, we can find a $\Lambda_x \in \C^*$
  such that $\Lambda_x (x) > 1$ and $\Lambda_x (\overline x) \leq 1$,
  for each $x \in \Min Q$.  Let
  $H_x \eqdef \Lambda_x^{-1} (]1, \infty])$: when $x$ varies over
  $\Min Q$, those open sets cover $\Min Q$, hence \blue{also its
    upward closure, which is $Q$}.  Since $Q$ is compact, there is a
  finite subset $E$ of $\Min Q$ such that
  $Q \subseteq \bigcup_{x \in E} H_x$.

  Let $\Lambda \eqdef \sum_{x \in E} \Lambda_x$, another element of
  $\C^*$.  Since $\varphi$ is linear by assumption,
  $\varphi (\Lambda) = \sum_{x \in E} \varphi (\Lambda_x)$.

  Relying on the compactness of $Q$,
  $\varphi (\Lambda) = \min_{z \in Q} \Lambda (z)$ is equal to
  $\Lambda (x_0)$ for some $x_0 \in Q$.  Since
  $Q \subseteq \bigcup_{x \in E} H_x$, there is an $x \in E$ such that
  $x_0 \in H_x$, namely such that $\Lambda_x (x_0) > 1$.  
  \BLUE{By the choice of $\Lambda_x$,}  $\Lambda_x (\overline x) \leq 1$, and
  therefore
  $\varphi (\Lambda_x) = \min_{z \in Q} \Lambda_x (z) \leq 1 <
  \Lambda_x (x_0)$.  For all the other points $y$ of $E$,
  $\varphi (\Lambda_y) = \min_{z \in Q} \Lambda_y (z) \leq \Lambda_y
  (x_0)$, so, summing up, we obtain that
  $\varphi (\Lambda) = \varphi (\Lambda_x) + \sum_{y \in E \diff
    \{x\}} \varphi (\Lambda_y) < \Lambda_x (x_0) + \sum_{y \in E \diff
    \{x\}} \Lambda_y (x_0) = \Lambda (x_0) = \varphi (\Lambda)$, a
  contradiction.
\end{proof}
\noindent 
Albeit natural, Lemma~\ref{lemma:jia} is surprising, for the following
reason.  \blue{We recall from Example~\ref{exa:cone:Val}~\ref{exa:C*}
  that $\C^*$ is given the weak$^*$upper topology, a.k.a.\ the
  weak$^*$-Scott topology.}  If one knows that $Q = \upc x$, then the
map $\varphi$ coincides with
\blue{$\text{ev}_x \colon \Lambda \mapsto \Lambda (x)$}, and must
therefore be continuous.  However, we are \emph{not} assuming that
$\varphi$ is continuous, only linear, in Lemma~\ref{lemma:jia}.  For
general non-empty convex compact saturated sets $Q$, we have no reason
to believe that
$\varphi \colon \Lambda \mapsto \min_{x \in Q} \Lambda (x)$ should be
lower semicontinuous from $\C^*$, with the weak$^*$-Scott topology, to
$\creal$.  But Lemma~\ref{lemma:jia} implies that it is, provided it
is linear.

\begin{thm}
  \label{thm:bary}
  Let $\C$ be a locally convex, locally convex-compact, sober
  topological cone $\C$.  Every continuous valuation $\nu$ on $\C$ has
  a unique barycenter $\beta (\nu)$, and the barycenter map $\beta$ is
  continuous; hence $\beta$ is the structure map of a
  $\Val_\wk$-algebra on $\C$, and the unique one that induces the
  given cone structure on $\C$.
\end{thm}
\begin{proof}
  Under the given assumptions on $\C$, $\Sm \C$ is a continuous d-cone,
  so Proposition~\ref{prop:cont:Valg} applies: there is a unique
  $\Val_\wk$-algebra structure
  $\alpha \colon \Val_\wk {(\Sm \C)} \to \Sm \C$ that induces the
  given cone structure on $\Sm \C$.

  Let us consider any continuous valuation $\nu$ on $\C$.  Since the
  unit $\eta^{\Sm}_\C$ of the $\Sm$ monad is, in particular,
  continuous, we can form the image valuation $\eta^{\Sm}_\C [\nu]$ on
  $\Sm \C$.  We let $Q \eqdef \alpha (\eta^{\Sm}_\C [\nu])$.

  For every $\Lambda \in \C^*$, by Lemma~\ref{lemma:Lambda:bar:Qcvx},
  $\overline\Lambda$ is in $(\Sm \C)^*$.  Since $\alpha$ is a
  barycenter map, we must have
  $\overline\Lambda (Q) = \int_{Q' \in \Sm \C} \overline\Lambda (Q')
  \dee\eta^{\Sm}_\C [\nu]$.  By the change of variable formula, the
  latter is equal to
  $\int_{x \in \C} \overline\Lambda (\eta^{\Sm}_\C (x))
  \dee\nu$\blue{, which is equal to}
  $\int_{x \in \C} \Lambda (x) \dee\nu$ \blue{by
    Lemma~\ref{lemma:Lambda:bar:Qcvx}}.  Also,
  $\overline\Lambda (Q) = \min_{x \in Q} \Lambda (x) = \varphi
  (\Lambda)$, taking the notation $\varphi$ from
  Lemma~\ref{lemma:jia}.  Since $\Lambda$ is arbitrary in $\C^*$, this
  shows that $\varphi$ is equal to the map
  $\Lambda \in \C^* \mapsto \int_{x \in \C} \Lambda (x) \dee\nu$,
  which is linear.  Hence, by Lemma~\ref{lemma:jia} (which applies
  since $\C$ is locally convex $T_0$ hence linearly separated; it is
  $T_0$ since sober), $Q = \upc x_\nu$ for some $x_\nu \in \C$.

  For every $\Lambda \in \C^*$, the equation
  $\overline\Lambda (Q) = \int_{x \in \C} \Lambda (x)\dee\nu$ then
  simplifies to
  $\Lambda (x_\nu) = \int_{x \in \C} \Lambda (x) \dee\nu$, showing
  that $x_\nu$ is a barycenter of $\nu$.  Since $\C$ is linearly
  separated and $T_0$, it is unique \cite{CEK:prob:scomp,GLJ:bary}.

  Let us define $\beta$ by $\beta (\nu) \eqdef x_\nu$, for every
  $\nu \in \Val_\wk \C$.  We claim that $\beta$ is continuous.  For
  every open subset $U$ of $\C$, for every $\nu \in \Val_\wk \C$,
  $\beta (\nu) \in U$ if and only if
  $\alpha (\eta^{\Sm}_\C [\nu]) \in \Box U$.  In other words,
  $\beta^{-1} (U) = (\alpha \circ \Val_\wk (\eta^{\Sm}_\C))^{-1} (\Box
  U)$, which is open since $\alpha$ and \blue{$\eta^{\mathcal S}_\C$} are continuous,
  and since $\Val_\wk$ is a functor.

  We conclude by Proposition~\ref{prop:Valg}.
\end{proof}

\blue{
  \begin{exa}
    \label{exa:thm:appl}
    We compare the strength of Theorem~\ref{thm:bary} with what we
    already know about barycenters in special cases (see
    Example~\ref{exa:bary}).
    \begin{enumerate}[label=(\arabic*)]
    \item For every core-compact space $X$, $\Lform X$ is a continuous
      d-cone, hence a sober, locally convex and locally
      convex-compact topological cone.  Therefore
      Theorem~\ref{thm:bary} applies; but \BLUE{Example~\ref{exa:bary}(1)
      suffices in this case.}
    \item Let $L$ be a sup-semilattice with the Scott topology and the
      cone structure defined in
      Example~\ref{example-semilattice-cones}~\ref{exa:cone:weird}. If
      $L$ is locally compact and sober in the Scott topology, then $L$
      is a sober, locally convex, locally convex-compact topological
      cone. Hence the barycenter map on $L$ exists, and the induced
      cone structure on $L$ is the original one on $L$. In this case,
      $L$ is actually a complete lattice, as sobriety of $L$
      guarantees that directed suprema on $L$ exist.  One can describe
      the barycenter map concretely, see
      Example~\ref{exa:bary}~\ref{exa:supp:weird}.
    \item Similarly with sup-semilattices with their upper topology
      (Example~\ref{exa:bary}~\ref{exa:supp:weird:upper}).
    \item\label{exa:C*:bary:yes} 
      A dual cone $\C^*$ is
      always topological, sober and locally convex (see
      Example~\ref{exa:cone:Val}~\ref{exa:C*}).
      Theorem~\ref{thm:bary} informs us that it suffices that $\C^*$
      be locally convex-compact to have a $\Val_\wk$-algebra structure
      inducing the cone structure on $\C^*$, and therefore for
      continuous valuations $\nu$ on $\C^*$ to have barycenters
      $\beta (\nu)$.  This barycenter is the map
      $x \mapsto \int_{\Lambda \in \C^*} \Lambda (x) \dee\nu$, see
      Example~\ref{exa:bary}~\ref{exa:C*:bary}. We will see in Example~\ref{dualcnobary}
      that not all continuous valuations $\nu$ on a dual cone $\C^*$ have  
      barycenters, hence not all $\C^*$ are locally convex-compact. 
    \item\label{exa:thm:appl:V} For every topological space $X$,
      $\Val_\wk X$ is topological, sober, and locally convex, see
      Example~\ref{exa:cone:Val}~\ref{exa:V}.  Hence
      Theorem~\ref{thm:bary} also gives us a $\Val_\wk$-algebra
      structure and barycenters if $\Val_\wk X$ is locally
      convex-compact.  We will see below that this holds true when $X$
      is stably locally compact, notably
      (Proposition~\ref{prop:Valg:locconvcomp}).  But the local
      convex-compactness assumption is unneeded, as we have seen in
      Example~\ref{exa:bary}~\ref{exa:bary:V}; examining this gap is
      the purpose of Section~\ref{sec:are-assumpt-necess}. 
    \end{enumerate}
\end{exa}}

\blue{%
\section{Are the assumptions necessary?}
\label{sec:are-assumpt-necess}
}
\noindent 
  One would like \blue{to know whether the assumptions of
  Theorem~\ref{thm:bary} are needed.  If a semitopological cone $\C$
  has a $\Val_w$-algebra structure that induces the given cone
  structure, then $\C$ has to be sober, topological, and weakly
  locally convex \cite[Proposition~4.9]{GLJ:bary}.  We do not know
  whether it has to be locally convex.  
  But does it have to be locally convex-compact?  
  Example~\ref{exa:thm:appl}~\ref{exa:thm:appl:V}
  suggests that the answer is no. However, we will give an example
  illustrating that Theorem~\ref{thm:bary} may fail if $\C$ is not 
  locally convex-compact.}
  
\blue{But first let us look at the special case
  of $\Val_\wk$-algebras in the category of stably compact spaces, as
  Cohen, Escard{\'o} and Keimel did \cite{CEK:prob:scomp}. In that setting, 
  we will show that it is indeed required that
  $\C$ be locally convex-compact.  We write $\SCcat$ for the category}
of stably compact spaces and continuous maps.  A space is \emph{stably
  locally compact} if and only if it is sober, locally compact, and
the intersection of any two compact saturated subsets is compact; it
is \emph{stably compact} if and only if it is also compact.  \blue{On
  a stably compact space $Y$, the complements of compact saturated
  subsets form a topology called the \emph{cocompact topology}.  The
  set $Y$ with the cocompact topology is then also stably compact
  \cite[Section~9.1.2]{JGL-topology}.}

\blue{In the proof of the following, we will use some results due to
  Gordon Plotkin \cite{Plotkin:alaoglu}.  Plotkin uses the
  weak$^*$-Scott topology on $\Val X$, equating it with $\C^*$ where
  $\C \eqdef \Lform X$; that is the same thing as the weak topology on
  $\Val X$, by \cite[Proposition~34]{AMJK:scs:prob}.  Let us write
  $[Q \geq r]$ for the collection of continuous valuations $\nu$ such
  that $\nu (U) \geq r$ for every open neighborhood $U$ of $Q$.  Then
  Corollary~1 of \cite{Plotkin:alaoglu}, together with subsequent
  comments, states that, for every stably locally compact space $X$,
  $\Val_\wk X$ is stably compact, and a subbase of closed subsets of
  the cocompact topology is given by the sets $[Q \geq r]$, where $Q$
  ranges over the compact saturated subsets of $X$ and $r \in \Rp$.
  In particular, all such sets $[Q \geq r]$ are compact in
  $\Val_\wk X$. They are clearly saturated.}


A \emph{linear retract} of a semitopological cone is a retract whose
retraction (not the section) is linear
\cite[Proposition~6.6]{heckmann96}.\\ 

\begin{prop}
  \label{prop:Valg:locconvcomp}
  For every stably locally compact space $X$, $\Val_\wk X$ is a
  locally convex-compact topological cone.  Each of those properties
  is preserved by linear retracts.  Hence, for every
  $\Val_\wk$-algebra $\alpha \colon \Val_\wk X \to X$ in $\SCcat$, $X$
  is a locally convex-compact topological cone in the induced cone
  structure.
\end{prop}
\noindent 
$X$ is also sober, since stably locally compact.  As a topological
cone, it will even be stably compact.  Hence only local convexity is
missing.

\begin{proof}
  We already know that $\Val_\wk X$ is a topological cone.  \blue{We
    use the following fact about locally compact spaces such as $X$:
    for every compact saturated subset $Q_0$ of $X$ and every open
    subset $U$ of $X$ \BLUE{that contains $Q_0$}, there is a compact saturated subset $Q$ of $X$
    such that $Q_0 \subseteq \interior Q \subseteq Q \subseteq U$
    \cite[Proposition~4.8.14]{JGL-topology}.  As a consequence, the
    family $\mathcal F_U$ of interiors $\interior Q$, where $Q$ ranges
    over the compact saturated subsets of $U$, is directed: it is
    nonempty by the previous observation with $Q_0 \eqdef \emptyset$,
    and for any two elements $\interior {Q_1}$ and $\interior {Q_2}$
    of $\mathcal F_U$, the same observation with
    $Q_0 \eqdef Q_1 \cup Q_2$ gives a larger element $\interior Q$ of
    $\mathcal F_U$.  Additionally, $U$ is the union of that family,
    since every point $x \in U$ is in some element of $\mathcal F_U$,
    by definition of local compactness.}
    
  Given $\nu \in \Val_\wk X$ and an open neighborhood $\mathcal U$ of
  $\nu$ \BLUE{in the weak topology}, $\mathcal U$ contains an open neighborhood of $\nu$ of the
  form $\bigcap_{i=1}^n [U_i > r_i]$, where each $U_i$ is open in $X$
  and $r_i > 0$.  Since $X$ is locally compact, \blue{each $U_i$ is
    the union of the directed family $\mathcal F_{U_i}$.  Since $\nu$
    is Scott-continuous, } $\nu (\interior {Q_i}) > r_i$ for some
  \blue{compact saturated subset $Q_i$ of $U_i$}.  Let $\epsilon > 0$
  be such that $\nu (\interior {Q_i}) > r_i+\epsilon$ for every $i$.
  Then $\nu$ is in $\bigcap_{i=1}^n [\interior {Q_i} > r_i+\epsilon]$.
  The latter is open and included in
  $\bigcap_{i=1}^n [Q_i \geq r_i+\epsilon]$, which is compact
  saturated\blue{, being a finite intersection of compact saturated
    sets in a stably compact space}.  Then,
  $\bigcap_{i=1}^n [Q_i \geq r_i+\epsilon]$ is included in
  $\bigcap_{i=1}^n [U_i > r_i] \subseteq \mathcal U$.  Finally,
  $\bigcap_{i=1}^n [Q_i \geq r_i+\epsilon]$ is convex, \blue{as an
    intersection of convex sets; any set of the form $[Q \geq r]$ is
    convex, since for all $\mu, \nu \in [Q \geq r]$, for every
    $a \in [0, 1]$, for every open neighborhood $U$ of $Q$,
    $(a \cdot \mu + (1-a) \cdot \nu) (U) = a \mu (U) + (1-a) \nu (U)
    \geq ar+(1-a)r = r$.} Therefore $\Val_\wk X$ is locally
  convex-compact.

  Let us consider a linear retraction $r \colon \C \to \Cb$, with
  associated section $s \colon \Cb \to \C$.  If $\C$ is locally
  convex-compact, then for every $y \in \Cb$ and every open
  neighborhood $V$ of \blue{$y$}, $s (y)$ has a convex compact neighborhood
  $Q$ included in $r^{-1} (V)$.  Let $U$ be the interior of $Q$.  The
  image of $Q$ under $r$ is convex, compact, contains the open
  neighborhood $s^{-1} (U)$ of $y$, and is included in $V$.  Hence
  $\Cb$ is locally convex-compact.  (This argument is the same as
  \cite[Proposition~6.6]{heckmann96}.)  Also, if $\C$ is topological,
  then addition and scalar multiplication on $\Cb$ are such that
  $x+y = r (s (x)+s(y))$ and $a\cdot x = r (a \cdot s (x))$, hence are
  jointly continuous, so $\Cb$ is topological.

  The final part \blue{follows} from the fact that any $\Val_\wk$-algebra is a
  linear retraction, with the unit $\eta_X$ as section.
\end{proof}
\noindent 
We now give an example to show that Theorem~\ref{thm:bary} may fail 
  if local convex-compactness of $\C$ is not assumed, even when $\C$ is a locally
  convex, sober topological cone.
  
 \begin{exa}\label{dualcnobary}
  Let $\real$ be the set of reals with the usual metric topology and $\C = \Val_\pw \real$,
  the set of point-continuous valuations on $\real$.  Similar to $\Val_\wk \real$,
  the space $\Val_\pw \real$, with the relative weak topology from $\Val_\wk \real$, also is a 
  locally convex sober topological cone~\cite{heckmann96}.
  
  The cone $\Val_\pw \real$ is actually strictly contained 
  in $\Val_\wk \real$. The Lebesgue valuation $\lambda$ on $\real$, obtained by restricting the 
  Lebesgue measure on opens of $\real$ is a continuous valuation that is not in $\Val_\pw \real$; see for example~\cite[Section 4.1]{heckmann96}. 
  The map $\eta_\pw \colon \real \to \Val_\pw \real \colon x\mapsto \delta_x$ is a well-defined continuous map, as $\delta_x$ 
  is point-continuous for each $x\in \real$. Hence $\eta_\pw[\lambda]$, the image valuation defined as $\Val_\wk(\eta_\pw)(\lambda)$, is a 
  continuous valuation in $\Val_\wk \Val_\pw \real$.  However, we claim that  $\eta_\pw[\lambda]$ does not have
  a barycenter in $\Val_\pw \real$. 
  
  Assume that $\eta_\pw[\lambda]$ does have a barycenter in $\Val_\pw \real$ and we denote it by~$\nu$. Then 
  by definition of barycenters, 
  $\Lambda (\nu) = \int_{\mu \in \Val_\pw \real} \Lambda(\mu)~d \eta_\pw[\lambda]$ for each continuous linear map 
  $\Lambda \colon \Val_\pw \real \to \overline \real_+$. In particular, each open subset $U$ in $\real$ determines 
  such a continuous linear map $\Lambda_U \colon \mu \mapsto \mu(U)$.  So for each open $U$ in $\real$, 
  \begin{align*}
   \nu(U) &= \Lambda_U (\nu) &\\
              &= \int_{\mu \in \Val_\pw \real} \Lambda_U(\mu)~d \eta_\pw[\lambda] &\\
              &= \int_{\mu \in \Val_\pw \real} \mu(U)~d \eta_\pw[\lambda] &\\
              &= \int_{x\in \real} \delta_x (U)~d \lambda & \text{(by the change of variable formula)}\\
              &= \int_{x\in  \real} \chi_U(x)~d \lambda & \\
               &=\lambda (U). &
   \end{align*}
  So $\nu = \lambda$. But this would contradict to the fact that $\lambda$ is not point-continuous, and hence $\eta_\pw[\lambda]$ does not 
  have any barycenters in $\Val_\pw \real$. As a byproduct, we know that  $\Val_\pw \real$ cannot be 
  locally convex-compact by Theorem~\ref{thm:bary}. 
  
  We also realize that $\Val_\pw \real$ is isomorphic to the dual cone $(\Lform \real)^*$, where $\Lform \real$ consists
  of all lower semi-continuous maps on $\real$, and both $\Lform \real$ and $(\Lform \real)^*$ are endowed with the 
  topology of pointwise convergence~\cite[Theorem 8.2]{heckmann96}. 
  Hence,  continuous valuations on a dual cone may fail to have barycenters. 
  \end{exa}

\section{Conclusion}
\label{sec:conclusion}

\blue{Theorem~\ref{thm:bary} is a pretty general barycenter existence
  theorem, which also serves to show that a large collection of
  topological cones has a (unique) $\Val_\wk$-algebra structure
  inducing its cone structure.  Through the examples we have given, we
  hope to have demonstrated that this theorem encompasses most of the
  cases where we already knew that barycenters existed, at the very least.}

We sum up the main remaining questions as follows.
\begin{itemize}
\item \blue{Which dual cones $\C^*$ are $\Val_\wk$-algebras?  We
    mentioned this in Example~\ref{exa:bary}~\ref{exa:C*:bary}, and we
    said that those that are locally convex-compact fit the bill in
    Example~\ref{exa:thm:appl}, and not all dual cones are locally
    convex-compact in Example~\ref{dualcnobary}.}
\item Given a $\Val_\wk$-algebra $\alpha \colon \Val_\wk X \to X$,
  \blue{we know that $X$ has an induced cone structure that is weakly
    locally convex.  Need it be locally convex?  Both local convexity
    and local convex-compactness are stronger properties than weak
    local convexity, and none is known to imply the other.  One should
    probably focus on the case of stably locally compact $X$, where we
    know at least that local convex-compactness is needed, by
    Proposition~\ref{prop:Valg:locconvcomp}.}
\item \blue{In the converse direction, what weaker conditions \BLUE{than}
    local convexity and local convex-compactness would be enough to
    deduce the same conclusions as those of 
  Theorem~\ref{thm:bary}?  Those conditions need to be at least as
  strong as weak local convexity.}
\end{itemize}


\section*{Acknowledgments}
\blue{Both authors would like to thank the anonymous referees for their careful reading of our manuscript and pointing out several areas where arguments were not clear. The final version has benefitted a lot from their comments.} Xiaodong Jia acknowledges support from NSFC (No.\;12371457, No.\;12231007). 

\bibliographystyle{alphaurl}
\bibliography{bary}

\end{document}